\documentclass[draft,11pt]{amsart}

\usepackage[top=2.5cm,bottom=2.5cm,left=2.5cm,right=2.5cm]{geometry}
\usepackage[english]{babel}
\usepackage{times}
\usepackage{dsfont}

\usepackage{amsmath,amssymb,hyperref,color,amsthm,verbatim}   
\usepackage[dvips]{graphicx}
\usepackage[mathscr]{euscript}
\usepackage{pb-diagram}
\usepackage[all]{xy}
\usepackage{cleveref}
\usepackage[colorinlistoftodos]{todonotes}
\usepackage{marginnote}

\newcommand{\R}{{\ensuremath{\mathbb{R}}}}

\newcommand{\ttimes}{\mathbin{\hbox to 0pt{$\times$\hskip -6.7pt\lower2.8pt\hbox{$\scriptscriptstyle\sim$}\hskip 0pt minus 1fil}\hphantom\times}}
\newcommand{\GLnR}{\textrm{GL}(\R^n)}
\newcommand{\framebundle}{\textrm{FR}(TM)}

\newcommand{\GL}{{{\rm{GL}}}}

\newcommand{\I}{{{\mathfrak{J}}}}

\newcommand{\Ad}{\mathrm{Ad}}
\newcommand{\g}{\mathfrak{g}}
\newcommand{\gl}{\mathfrak{gl(\R^n)}}

\title[Infinitesimally homogeneous manifolds]{Infinitesinally homogeneous manifolds with prescribed structure groups}

\author[C. Mar\'in]{Carlos Alberto Mar\'in Arango}
\author[D. Bl\'azquez-Sanz]{David Bl\'azquez-Sanz}

\address{Instituto de Matem\'aticas \hfill\break\indent  Universidad de Antioquia \hfill\break\indent Medell\'in, Colombia}
\email{calberto.marin@matematicas.udea.edu.co}
\address{Universidad Nacional de Colombia - Sede Medell\'in \hfill\break\indent  Facultad de Ciencias 
\hfill\break\indent Escuela de Matem\'aticas \hfill\break\indent  Medell\'in, Colombia}
\email{dblazquezs@unal.edu.co}

\subjclass[2000]{53A15, 53B05, 53C10, 53C30}

\keywords{Infinitesimally homogeneous manifold, Inner torsion, $\mathbf G$-structures}

\date{November 2015}

\begin{document}

\makeatletter
\newenvironment{dem}{Proof:}{\qed}

\makeatother

\numberwithin{equation}{section}
\theoremstyle{plain}\newtheorem{teo}{Theorem}[section]
\theoremstyle{plain}\newtheorem{prop}[teo]{Proposition}
\theoremstyle{plain}\newtheorem{lema}[teo]{Lemma}

\theoremstyle{definition}\newtheorem{df}[teo]{Definition}
\theoremstyle{remark}\newtheorem{rem}[teo]{Remark}

\theoremstyle{definition}\newtheorem{example}[teo]{Example}

\begin{abstract}
We explore the class of triples $(M,\nabla,P)$ where $M$ is a manifold, $\nabla$ is an affine connection in $M$ and $P$ is a $G$-structure in $M$. Inside this class there are infinitesimally homogeneous manifolds, characterized by having $G$-constant curvature, torsion and inner torsion. For each matrix Lie group $G\subseteq {\rm GL}(\mathbb R^n)$ there is a class of infinitesimally homogeneous manifolds with structure group $G$. In this paper we characterize the classes of infinitesimally homogeneous manifolds for some specific values of the structure group $G$ including: identity group, finite groups, diagonal group, special linear group, orthogonal group and unitary group.
\end{abstract}

\maketitle

\section{Introduction}
Let $G\subset \mathrm{GL}(\mathbb{R}^n)$ be a Lie subgroup. By a \emph{manifold with affine connection and $G$-structure} we mean a triple $(M,\nabla,P)$ in which $M$ is an $n$-dimensional manifold, $\nabla$ is an affine connection on $M$ and $P\subset\framebundle$ is a $G$-structure on $M$. Throughout this paper, by an affine connection on $M$ we mean a linear connection on the tangent bundle $TM$. The geometry of a manifold with affine connection and $G$-structure $(M,\nabla,P)$ is described by three tensors: the curvature tensor $(R)$ of $\nabla$, the torsion tensor $(T)$ of $\nabla$ and the inner torsion $(\mathfrak{J}^P)$ of $P$, see \cite{tauskpaper}. 

In the case of a manifold endowed with a $G$-structure $P$, the notion of \emph{$P$-constant tensor field} makes sense, and it refers to a tensor field whose representation in frames that belong to $P$ is constant; in the sense of being independent of the particular frame and the point on the manifold. 

A manifold with affine connection and $G$-structure $(M,\nabla,P)$ is said to be \emph{infinitesimally homogeneous} if $R$, $T$ and $\mathfrak{J}^P$ are all $P$-constant. That is, there exist multilinear maps $R_0:\R^n\times\R^n\times\R^n\to\R^n,$ $T_0:\R^n\times\R^n\to\R^n$ and a linear map $\mathfrak{J}_0:\R^n\to\gl/\g,$ $\g$ being the Lie algebra of $G$ such that
\begin{equation}\label{eq:homogeneidadinfinitesimal}
\begin{gathered}
p^*R_x=R_0,\quad p^*T_x=T_0,\\
\overline{\mathrm{Ad}_p}\circ\mathfrak{J}_0=\mathfrak{J}^P_x\circ p,
\end{gathered}
\end{equation}
for each $x\in M$, and each $p\in P_x$. The maps $R_0$, $T_0$ and $\mathfrak{J}_0$ are collectively called the \emph{characteristic tensors} of $(M,\nabla,P)$, because they provide a local characterization of infinitesimally homogeneous manifolds in the sense that two of these having the same characteristic tensors are locally equivalent by means of a connection and $G$-structure preserving diffeomorphism. Infinitesimally homogeneous manifolds have been studied thoroughly by Piccione and Tausk in \cite{tauskpaper}, in the aforementioned paper the authors proved an existence result for local and global immersions into infinitesimally homogeneous manifolds. This is a very general result and includes several isometric immersions theorems that appear in the literature. From now on, by \emph{infinitesimally homogeneous manifold} we mean infinitesimally homogeneous manifold with affine connection and $G$-structure.


In \cite{marin}, the first author gave necessary and sufficient conditions for maps $R_0$, $T_0$ and $\mathfrak{J}_0$ to be the characteristic tensors of an infinitesimally homogeneous manifold $(M,\nabla, P)$. More specifically, given a Lie subgroup $G\subset \GLnR$, the main result of \cite{marin} is an algebraic characterization (involving both the group $G$ and its Lie algebra $\mathfrak{g}$) for the possible characteristic tensors of an infinitesimally homogeneous  manifold $(M,\nabla, P)$ with structure group $G$. Therefore, for a given Lie subgroup $G\subset \GLnR$, the
problem of classifying the infinitesimally homogeneous manifolds with structure group $G$ is reduced to the problem of classifying all the maps $R_0, T_0$ and $\mathfrak{J}_0$ that satisfy the algebraic conditions given in \cite{marin}. As an application, in this paper we give specific characterizations of infinitesimally homogeneous manifolds for certain Lie groups.

\section{Notation and Preliminaries}

\subsection{Inner torsion}
Let $G\subset \GLnR$ be a Lie subgroup whose Lie algebra will be denoted by $\mathfrak{g}$, and let $(M,\nabla, P)$ be a manifold with affine connection and $G$-structure. If $\omega$ denotes the $\mathfrak{gl}(\R^n)$-valued connection form on $\mathrm{FR}(TM)$ associated with $\nabla$, it is clear that $\nabla$ is {\em{compatible\/}} with the subbundle $P$ if $\omega|_P$ is $\mathfrak{g}$-valued. In order to handle the general case in which $\nabla$ is not compatible with $P$, the concept of {\em{inner torsion\/}} was introduced in \cite{tauskpaper}. It is a tensor that measures the lack of compatibility of $\nabla$ with $P$. This tensor may be defined in a briefly way as follows: for each $x\in M$, let $G_x$ be the subgroup of $\GL(T_xM)$ consisting of all {\em $G$-structure preserving maps}, i.e., maps $\sigma:T_xM\to T_xM$ such that $\sigma \circ p\in P_x$ for some $p\in P_x$. Clearly, $G_x= \mathcal{I}_p(G)$, for each $p \in P_x$ where $\mathcal{I}_p:\GL(\R^n) \to \GL(T_xM)$ denotes the Lie groups isomorphism given by conjugation by $p$, so that $G_x$ is a Lie subgroup of $\GL(T_xM)$ whose Lie algebra we denote by $\mathfrak{g}_x$. Note that the linear isomorphism $\mathrm{d}\mathcal{I}_p(1)= \mathrm{Ad}_p$ carries $\mathfrak{g}$ onto $\mathfrak{g}_x$, so it induces an isomorphism of the quotients $\overline{\mathrm{Ad}}_p:\mathfrak{gl}(\R^n)/\mathfrak{g}\to \mathfrak{gl}(T_xM)/\mathfrak{g}_x$. Let $s:U\subset M\to P$ be a smooth local section of $P$ around a point $x$ and set $\overline{\omega}=s^*\omega$, $s(x)=p$. The map
\begin{equation}\label{eq:e2}\xymatrix{%
T_xM \ar@.@/_1.5pc/[rrr]_{\mathfrak{I} ^P_x}\ar[r]^-{\overline{\omega}_x}&
\mathfrak{gl(\R^n)} \ar[r]^-{\mathfrak q}& \mathfrak{gl(\R^n)}/\mathfrak{g}
\ar[r]^-{\overline{{\rm{Ad}}}_{p}}& \mathfrak{gl}(T_xM)/\mathfrak g_x }
\end{equation}
does not depend on the choice of the local section $s$ (see \cite{tauskpaper}). The linear map $\mathfrak{I}^P_x$ defined by
\eqref{eq:e2} is called the {\em inner torsion\/} of the $G$-structure $P$ at the point $x\in M$  with respect to the connection $\nabla$.
If $\Gamma$ denotes the Christoffel tensor of the connection $\nabla$ with respect to the local frame $s$, then $\Gamma _x =\Ad _p \circ \overline{\omega}_x$, where $p=s(x)$. Thus the inner torsion $\mathfrak{I} ^P_x:T_xM \to \mathfrak{gl}(E_x)/\mathfrak g_x$ of the $G$-structure $P$ at the point $x$ is given by the composition of the Christoffel tensor $\Gamma_x:T_xM\to \mathfrak{gl}(T_xM)$ of $\nabla$ with respect to $s$ and the quotient map $\mathfrak{gl}(T_xM)\to \mathfrak{gl}(T_xM)/\mathfrak g_x$.

\subsection{Invariance properties of characteristic tensors}\label{sec:condicionesalgebraicas}
Let $(M,\nabla,P)$ be a manifold with affine connection and $G$-structure and let $\mathfrak g$ be the Lie algebra of $G$. Let $T_0$, $R_0$, $\mathfrak{J}_0$ be its characteristic tensors. By definition, $T_0$, $R_0$, $\mathfrak{J}_0$ are invariant for the natural action of $G$. That is, they satisfy the following relations:
\begin{gather}
\label{eq:R0}
R_0(u,v)=\Ad_g\cdot R_0(g^{-1}\cdot u,g^{-1}\cdot v);\\
\label{eq:T0}
T_0(u,v)=g\cdot T_0(g^{-1}\cdot u,g^{-1}\cdot v);\\
\label{eq:I0}
\Ad_g\big(\lambda(u)\big)-\lambda(g\cdot u)\in \g,
\end{gather} 
for all $g\in G$ and $u, v\in \R^n$. Here $\lambda: \R^n\to\mathfrak{gl}(\R^n)$ is an arbitrary lift of $\I_0$. Taking
$g = \exp(tL)$ and differentiating the above expressions with respect to $t$ we obtain the following relations:
\begin{gather}
\label{eq:R02}
[L, R_0(u,v)]-R_0(L\cdot u,v)-R_0(u,L\cdot v) =0;\\
\label{eq:T02}
L\circ T_0(u,v)-T_0(L\cdot u,v)-T_0(u,L\cdot v)=0;\\
\label{eq:I02}
[L,\lambda(u)]-\lambda(L\cdot u)\in \g.
\end{gather}
that hold for each $L\in \mathfrak{g}$ and $u, v \in \R^n$. 
If $G$ is a connected Lie group then relations \cref{eq:R0,eq:T0,eq:I0} and
\cref{eq:R02,eq:T02,eq:I02} are mutually equivalent.

Algebraic necessary and sufficient conditions for multilinear maps $T_0$, $R_0$ and $\mathfrak{J}_0$ to be the characteristic tensors of an infinitesimally homogeneous manifold, as well as some additional relations and properties, are given in \cite[Theorems 1 and 2]{marin}.





\section{Infinitesimally homogeneous manifolds with prescribed structure groups}

Here we explore the classes of infinitesimally homogeneous manifolds with specific structure groups. This means that for a given Lie group $G\subset \mathrm{GL}(\R^n)$ we find the possible characteristic tensors $R_0$, $T_0$ and $\mathfrak{J}_0$ of infinitesimally homogeneous manifolds with structure group $G$. This allows us to characterize the class of infinitesimally homogeneous manifolds with structure group $G$ is terms of the language of classical differential geometry. 
From the facts mentioned in section \ref{sec:condicionesalgebraicas} we know that a first step is to find all the 
possible multilinear maps $T_0, R_0, \mathfrak{J}_0$ satisfying \cref{eq:R0,eq:T0,eq:I0}.

\begin{lema}\label{thm:lemasbasicos}
Let $(M,\nabla,P)$ be an infinitesimally homogeneous manifold with characteristic tensors $T_0$, $R_0$,
$\mathfrak J_0$ and structure group $G \subset {\rm GL}(\mathbb R^n)$. The following statements hold:
\begin{enumerate}
\item
If $G$ contains more than one scalar matrices, then the multilinear map $T_0$ 
necessarily vanish.
\item
If $G$ contains more than two scalar matrices, then the multilinear maps $R_0$ and $T_0$ 
necessarily vanish. 
\item
If $-1\in G$ then the multilinear maps $T_0$ and $\mathfrak{J}_0$ 
necessarily vanish. 
\end{enumerate} 
\end{lema}
\begin{dem}
Everything follows directly from relations \eqref{eq:R0}, \eqref{eq:T0} and \eqref{eq:I0}. In order to illustrate this idea we give a complete proof of (3). Taking $g=-1$ in \eqref{eq:T0} we obtain $T_0(u,v)=-T_0(u,v)$ and thus $T_0=0$. Also, if $\lambda:\R^n\to \mathfrak{gl}(\R^n)$ is an arbitrary lift of $\mathfrak{J}_0$ then from \eqref{eq:I0} it follows:
\[
\lambda\left(g(u)\right) - g\circ \lambda (u) \circ g^{-1} \in \mathfrak{g}
\] 
for all $g\in G$ and $u\in \R^n$. In particular, for $g=-1$ we obtain  $-2\lambda (u) \in \mathfrak{g}$ for any $u\in \R^n$. That
is, $\mathfrak{J}_0 =0$.
\end{dem}

\subsection{Trivial structures.} Let us consider the case $G=\mathrm{GL}(\R^n)$. Here, there is no choice for the $G$-structure, it is the whole frame bundle $\mathrm{FR}(TM)$. Let $(M,\nabla,\mathrm{FR}(TM))$ be an infinitesimally homogeneous manifold. Its inner torsion automatically vanish. Moreover, there are no non zero $G$-invariant tensors $T_0$ and $R_0$ so that $\nabla$ is flat and torsion free. Finally, $(M,\nabla,\mathrm{FR}(TM))$ is an infinitesimally homogeneous manifold if and only if $(M,\nabla)$ is an affine manifold in the sense of Thurston (see \cite{HT}).

\subsection{Infinitesimally homogeneous structures over the identity.}

Let us consider the case $\{1\}\subset \mathrm{GL}(\R^n)$. By definition, a $\{1\}$-structure in $M$ is
a global frame $P = \{X_1,\ldots,X_n\}$. Let $\nabla$ be a linear connection in $M$, and let us assume that $(M,\nabla,P)$ is an infinitesimally homogeneous manifold with structure group $\{1\}$. 
The inner torsion $\mathfrak{J}_{x}^P :T_xM\to \mathfrak{gl}(T_xM)$ 
coincides with the Christoffel tensor $\Gamma_x:T_xM\to \mathfrak{gl}(T_xM)$ of the connection $\nabla$ with respect to 
the frame $P_x$. Let us consider $\{\theta_1,\ldots,\theta_n\}$ the dual co-frame associated to $P$. The 
expression in coordinates of the inner torsion is:
$$\mathfrak{J}^P = \Gamma_{ij}^k \theta_i \otimes \theta_j \otimes X_k,$$ 
where the functions $\Gamma_{ij}^k$ are the so called Christoffel symbols of $\nabla$ in the frame $P$. The infinitesimal homogeneity implies that the $\Gamma_{ij}^k$ are constant. 

Since $P$ is a global frame, we can find structure functions $\lambda_{ij}^k$ such that,
$[X_i,X_j]=\lambda_{ij}^kX_k$. The expression in coordinates of the torsion $T$ of $\nabla$ yields:

$$T = (\Gamma_{ij}^k - \Gamma_{ji}^k - \lambda_{ij}^k)\theta_i\wedge \theta_j \otimes X_k.$$
As before, the infinitesimal homogeneity implies that the components of $T$ are constant functions, and thus
the $\lambda_{ij}^k$ are also constant. Thus, the frame $P$ is a \emph{paralellism of $M$} modeled over a $n$-dimensional
Lie algebra. It is also clear that if Christoffel symbols $\Gamma_{ij}^k$ and the structure functions $\lambda_{ij}^k$ 
are constants, then the curvature tensor $R$ is also $P$-constant. Let us denote by $\mathfrak h$ the Lie
algebra of vector fields spanned by $X_1,\ldots,X_n$. The Christoffel tensor can be seen
as a linear map from $\mathfrak h$ to $\mathrm{End}(\mathfrak h)$. Summarizing:

\begin{prop}
A infinitesimally homogeneous manifold $(M,\nabla,P)$ with structure group $\{1\}$ is characterized by the following
data on the manifold $M$:
\begin{enumerate}
\item A global frame $P = \{X_1,\ldots,X_n\}$ such that the vector fields $X_1,\ldots,X_n$ span a $n$-dimensional Lie algebra
$\mathfrak h$ of vector fields in $M$.
\item A linear map $\Gamma \colon \mathfrak h \to \mathrm{End}(\mathfrak h)$.
\end{enumerate}
\end{prop}

The standard model for an infinitesimally homogeneous manifold with structure group $\{1\}$
is a Lie group $H$, endowed of a basis of its
Lie algebra $\mathfrak h$ of left invariant vector fields and an endomorphism from $\mathfrak h$ to 
$\mathrm{End}(\mathfrak h)$.
In particular, if $H$ is semisimple and we consider the adjoint action as such linear map, we obtain the Levi-Civita connection of the Killing metric in the Lie group. Also, by the third Lie theorem we know that any $n$-dimensional manifold with a transitive action of an $n$-dimensional Lie algebra is locally isomorphic to a Lie group with the action of its Lie algebra. Thus, any infinitesimally homogeneous affine manifold with structure group $\{1\}$ is locally isomorphic to a Lie group $H$, with a given basis of its Lie algebra $\mathfrak h$ and a given linear map from 
$\mathfrak h$ to $\mathrm{End}(\mathfrak h)$.

\subsection{Infinitesimally homogeneous structures over finite groups}
Let us consider $G \subset \mathrm{GL}(\R^n)$ a finite matrix group and
$(M,\nabla, P)$ an infinitesimally homogeneous manifold with structure group $G$. 
By definiton, the projection $\pi\colon P \to M$
is a covering space. The group $G$ acts freely
and transitively in each fiber of $\pi$, the quotient $P/G$ is identified with 
$M$, and $\pi$ is a Galois cover space modeled over the finite group $G$. 
Let us consider $\bar P$ the fibered product of coverings $P\times_M P$, endowed with the projection $\overline\pi:\overline P\to  P$ defined by
$\bar\pi(p_x,q_x) = p_x$. By definition of Galois cover, the projection 
$\bar\pi\colon \bar P \to P$ is a trivial cover. The diagonal set:
$$P' = \{(p_x, p_x) \colon p_x\in P\} \subset \bar P$$
is a connected component of $\bar P$. The rest of connected components of $\bar P$
are obtained by the action of $G$ on the frames:
$$P' g = \{(p_x, p_xg)   \colon p_x\in P\}.$$
Since $\pi$ is a local diffeomorphism, for each $p_x\in P$ it induces an isomorphism between 
$\textrm{FR}(T_{p_x}P)$ and $\textrm{FR}(T_xM)$. By means
of these isomorphisms we have that $P'$ can be seen as global frame
and then a $\{1\}$-structure in $P$. Same facts happen for $P'g$ for each $g\in G$
and the triples $(P,\pi^*(\nabla),P'g)$ are a family of
infinitesimally homogeneous structures in $P$ over $\{1\}$.
The group $G$ acts in $P$ by diffeomorphisms. Let us write $\phi_g(p_x)$ for $p_xg$. The
diffeomorphisms $\phi_g$ are symmetries of the connection $\pi^*(\nabla)$. We may
differentiate $\phi_g$ to compute the action of $G$ in the global frames, and we obtain
$(\phi_g)_*(P') = P'g^{-1}$. The triple $(P,\pi^*(\nabla),P')$ together with the action
of $G$ in $P$ completely describes the 
infinitesimally homogeneous manifold $(M,\nabla,P)$ as the quotient of $P$ by $G$.
We may then state the following result that characterizes the infinitesimally homogeneous
manifolds with structure group $G$ as quotients of infinitesimally homogeneous manifolds with structure
group $\{1\}$ by suitable groups of symmetries of their respective affine connections.

\begin{prop}
Let us consider $G \subset \mathrm{GL}(\R^n)$ a finite matrix group
and $(M,\nabla,P)$ an infinitesimally homogeneous manifold with structure group $\{1\}$
endowed with right action of $G$ in $M$ satisfying:
\begin{enumerate}
\item $G$ acts in $M$ freely and completely discontinuously by symmetries of $\nabla$.
\item $G$ acts equivariantly in the global frame $P$ in the following sense, for
each $g\in G$, $(\phi_g)_*(P) = Pg^{-1}$.
\end{enumerate}
The quotient map $\pi\colon M\to M/G$ is a cover and there is a projected
connection $\pi_*(\nabla)$ in $M$. The triple $(M/G,\pi_*(\nabla),P)$ is an infinitesimally homogeneous manifold
with structure group $G$. Conversely, any infinitesimally homogeneous manifold with structure group $G$ is isomorphic to a quotient manifold of this kind. 
\end{prop}

\subsection{Oriented Riemannian structures}
Let us consider the group $G={\rm SO}(n)$ consisting of all orthogonal matrices whose determinant is $1$. For this group a $G$-structure on a smooth manifold $M$ can be thought as having an orientation of $M$ and a Riemannian structure. 
In order to determine the suitable candidates for the characteristic tensors $R_0$, $T_0,$ and $\mathfrak{J}_0$ in what follows we will show some algebraic results.

\begin{lema}\label{curvature1}
Let $R$ be a quadrilinear map $R\colon \mathbb R^n \times \mathbb R^n \times \mathbb R^n \times \mathbb R^n \longmapsto \mathbb R$ satisfying:
\begin{itemize}
\item[(a)] $R(u_1,u_2,u_3,u_4) = -R(u_1,u_2,u_3,u_4),$
\item[(b)] $R(u_1,u_2,u_3,u_4) = -R(u_1,u_2,u_4,u_3),$
\item[(c)] $R(u_1,u_2,u_3,u_4) + R(u_1,u_3,u_4,u_2) + R(u_1,u_4,u_2,u_3) = 0,$
\end{itemize}
for all $u_1,u_2,u_3,u_4\in \mathbb R^n$. If $R$ is ${\rm SO}(n)$-invariant then $R$ is a multiple of the tensor $K$,
\begin{equation}\label{eq:curvaturaconstante} K:(u_1,u_2,u_3,u_4)\longmapsto \langle u_2,u_3\rangle \langle u_1,u_4\rangle -\langle u_1,u_3 \rangle \langle u_2,u_4 \rangle.\end{equation} 
\end{lema}
\begin{dem}
Let us take $\lambda = R(e_1,e_2,e_1,e_2)$. For any orthogonal basis $\{u,v\}$ of a plane in $\mathbb R^n$ there is a rotation $g$ such that $g(e_1) = u$, $g(e_2) = v$ (we may change the order of $u,v$ in the case $n=3$). By the ${\rm SO}(n)$-invariance of $R$ we have:
$$R(u,v,u,v) = R(e_1,e_2,e_1,e_2) = \lambda = \lambda K(u,v,u,v).$$
Finally, by Lemma \cite[Chapter 5, Proposition 1.2]{KN} we have that $R = \lambda K$.
\end{dem}

\begin{rem}
In the case $n\neq 4$ it is possible to proof that the only ${\rm SO}(n)$-invariant $4$-linear maps satisfying condition (a) are the multiples of $K$. 
\end{rem}

It is clear that, by contraction with the scalar product, the ${\rm SO}(n)$-invariant $4$-linear maps are in 1-1 correspondence with the  ${\rm SO}(n)$-equivariant $(3,1)$-tensors. Therefore, the characteristic tensor $R_0$ of an infinitesimally homogeneous manifold with structure group ${\rm SO}(n)$ is an scalar multiple of the $(3,1)$-tensor $K_0$ corresponding to $K$, which defined by the formula $K_0(u,v,w) = \langle v,w \rangle u - \langle u, w \rangle v$. Note that $\lambda K_0$ is the curvature tensor of a Riemannian manifold of constant sectional curvature equal to $\lambda$.
\smallskip

For the characteristic tensor $T_0$ we have the following algebraic result. 

\begin{lema}\label{lema36}
For $n\neq 1,3$ the space of skew-symmetric ${\rm SO}(n)$-invariant $(2,1)$-tensors in $\mathbb R^n$ vanishes. For $n = 3$ this space is spanned by the euclidean vector product:
$$\wedge \colon \mathbb R^3 \times \mathbb R^3 \to \mathbb R^3, \quad (u,v) \mapsto u\wedge v$$
\end{lema}

\begin{dem}
The case $n=3$ corresponds to the uniqueness of the vector product in euclidean space, which is already well known. 
If $n$ is even, then $-1\in {\rm SO}(n)$ and the result follows from Lemma \ref{thm:lemasbasicos}. On the other hand, if $n$ is odd and $n\ge 5$, let us consider an orthonormal basis $B=\{b_1,\dots,b_{n}\}$ of $\R^n$. The $G$-invariance of $T_0$ is equivalent to the identity
\begin{equation}\label{eq:torsionorientada}
g\left ( T_0(u,v)\right) = T_0\left(g(u),g(v)\right),
\end{equation} for all $g\in G$ and $u,v \in \mathbb R^n.$
Given $b_k, b_l \in B$, let $g_{kl}\in {\rm SO}(n)$ be the element given by the rotation of angle $\pi/2$ in the $(b_k,b_l)$-plane. 
Let us consider fixed indexes $i, j \in \{1,\dots, n\}$.
Replacing $g$ by $g_{kl}$ in the identity \eqref{eq:torsionorientada}, we get,
\[
g_{kl}\big(T_0(b_i,b_j)\big)=T_0\big(g_{kl}(b_i),g_{kl}(b_j)\big),\] for all $k,l \in \{1,\dots , n\}$. 

\noindent In the case if repeated indexes $i=k, j=l$ we obtain
\[
g_{ij}\big(T_0(b_i,b_j)\big)=T_0\big(g_{ij}(b_i),g_{ij}(b_j)\big)= T_0\big(b_i,b_j\big),
\]consequently $T_0(b_i,b_j)\in \mathrm{gen}(b_i,b_i)^{\perp}$.

\noindent For $k, l \notin \{i,j\}$, 
\[
g_{kl}\big(T_0(b_i,b_j)\big)=T_0\big(g_{kl}(b_i),g_{kl}(b_j)\big)=T_0\big(b_i,b_j\big),\] so that $T_0(b_i,b_j)\in \mathrm{gen}(b_k,b_l)^{\perp}$. Hence, $\langle T_0(b_i,b_j),b_k\rangle=0$ for all $k\in \{1,\dots, n\}$ Since $b_i, b_j$ are arbitrary elements in $B$, the result follows. 
\end{dem}

Finally let us proceed to determine the suitable candidates for the characteristic tensor $\mathfrak J_0$.

\begin{lema}\label{lema37}
Let us consider the vector space $\mathfrak{gl}(n)/\mathfrak{so}(n)$ endowed with the adjoint action of $\rm{SO}(n)$. Then, the space of linear $SO(n)$-invarivariant maps from $\mathbb R^n$ to $\mathfrak{gl}(n)/\mathfrak{so}(n)$ vanishes. In other words the tensor $\mathfrak J_0$ vanishes. 

\end{lema}
\begin{dem}
If $n$ is even, then $-1\in {\rm SO}(n)$ and the result follows from Lemma \ref{thm:lemasbasicos}. On the other hand, we assume that $n$ is odd and $n\ge 3$. Let $\mathfrak{I}$ be a $SO(n)$-invariant linear map from $\mathbb R^n$ to $\mathfrak{gl}(n)/\mathfrak{so}(n)$, and let $\lambda:\R^n \to \mathfrak{gl}(\R^n)$ be an arbitrary lifting for $\mathfrak{I}$. In order to obtain the desired result we will show that $\lambda$ is $\mathfrak{so}(n)$-valued.
To do this, let $B=\{b_1,\dots, b_n\}$ be an orthonormal basis of $\mathbb R^n$ and let us consider a fixed vector $b_i\in B$. For each $g\in {\rm SO}(n)$, we will write $L_g$ to denote the linear map  
\[L_g:=\Ad_{g}\big(\lambda (b_i)\big)-\lambda\big(g(b_i)\big) \in \mathfrak{so}(n).
\]  
The $SO(n)$-invariance \eqref{eq:I0} implies that $L_g \in \mathfrak{so}(n)$ for all $g\in {\rm SO}(n)$. Therefore for each $j,k \in \{1,\dots, n\}$ and each $g\in {\rm SO}(n)$ we have that
\begin{equation}\label{torsioninternaorientadaymetrica}
\langle L_g (b_j),b_k\rangle = - \langle L_g (b_k),b_j\rangle,
\end{equation} for all $g\in {\rm SO}(n)$. In particular, 
$\langle L_g(b_i),b_i\rangle =0,$ for all $g\in {\rm SO}(n)$. 
Thus given $g\in {\rm SO}(n)$ such that $g(b_i)=-b_i$, we have 
\[
\langle g\left(\lambda(b_i)\right)\cdot (-b_i),b_i\rangle = - \langle \lambda(b_i)\cdot b_i,b_i\rangle;
\]
consequently,
\begin{equation}\label{eq:torsioninternametricaorientadanula1}
\langle \lambda(b_i)\cdot b_i,b_i\rangle =0.
\end{equation}
On the other hand given $b_j, b_k \in B$, let us denote by $g_{jk}$ the element in ${\rm SO}(n)$ obtained by rotation of angle $\pi/2$ in the $(b_j,b_k)$-plane. Replacing $g$ by $g_{jk}$ in \eqref{torsioninternaorientadaymetrica} results
\[
\langle \lambda(b_i)\cdot b_j,b_k\rangle + \langle \lambda\left(g_{jk}(b_i)\right)\cdot b_j,b_k\rangle = -\langle \lambda(b_i)\cdot b_k,b_j\rangle - \langle \lambda\left(g_{jk}(b_i)\right)\cdot b_k,b_j\rangle;
\] for $i\ne j,k$,
\begin{equation}\label{eq:torsioninternametricaorientadanula2}
\langle\lambda(b_i)\cdot b_j,b_k\rangle=-\langle \lambda(b_i)\cdot b_k,b_j\rangle.
\end{equation}
In particular we get
\begin{equation}\label{eq:torsioninternametricaorientadanula3}
\langle \lambda(b_i)\cdot b_j,b_j\rangle =0,
\end{equation} for all $j\ne i$.
By combining \eqref{eq:torsioninternametricaorientadanula1} and \eqref{eq:torsioninternametricaorientadanula3} we obtain 
\[
\langle \lambda(b_i)\cdot (b_i+b_j),b_i+b_j\rangle =0,
\]for all $j\in \{1,\dots, n\}$, which implies that
\begin{equation}\label{eq:torsioninternametricaorientadanula4}
\langle \lambda(b_i)\cdot b_j,b_i\rangle = - \langle \lambda(b_i)\cdot b_i,b_j\rangle;
\end{equation} for all $j\ne i$.
From equalities \eqref{eq:torsioninternametricaorientadanula2} and \eqref{eq:torsioninternametricaorientadanula4} we clonclude that $\lambda (b_i)\in \mathfrak{so}(n)$ and the result follows.
\end{dem}

\begin{teo}
Let $(M,\nabla,P)$ be an infinitesimally homogeneous manifold of dimension $n\neq 1$ with structure group $\mathrm{SO}(n)$ . Then $P$ is the bundle of positively oriented orthogonal frames in an oriented Riemannian structure in $M$ with constant sectional curvature. 
Additionally, one of the following cases holds:
\begin{itemize}
\item[(a)]  $\nabla$ is the Levi-Civita connection of the Riemannian structure in $M$.
\item[(b)] The dimension of $M$ is $n=3$, and $\nabla$ is of the form $\nabla'+\mu\mathfrak{t}$ were $\nabla'$ is the Levi-Civita connection, $\mu$ is an scalar constant, and $\mathfrak{t}$ is the vector product in $TM$ induced by the oriented Riemannian structure of $M$. 
\end{itemize}
\end{teo}

\begin{dem}
In any case, by Lemma \ref{curvature1}, the characteristic tensor $R_0$ is of the form $\lambda K_0$ and thus $(M,P)$ is 
a Riemannian oriented manifold of constant sectional curvature. Let us now consider the case $n\neq 3$. In 
such case, by Lemmas \ref{lema36} and \ref{lema37}, the characteristic tensors $T_0$, $\mathfrak J_0$ vanish. 
It follows that $\nabla$ is
a torsion-free connection parallel to the oriented Riemannian structure, and thus it is its Levi-Civita connection. Finally let us consider the case $n=3$. By Lemma \ref{lema37} we have that $\nabla$ is a connection parallel to the oriented Riemann
structure. Let us consider $\nabla'$ the torsion-free part of $\nabla$ defined by $\nabla' = \nabla - \frac{1}{2}T$.
By the above argument, $\nabla'$ is the Levi-Civita connection in the Riemannian structure, and finally, by Lemma
3.6 the torsion $T$ is a multiple of the vector product $\mathfrak t$ induced by the oriented Riemannian structure $P$
in $TM$.
\end{dem}

\begin{rem}
Note that ${\rm SO}(1) = \{1\}$. Therefore, the case $n=1$ has been already discussed as an infinitesimally homogeneous structure over the identity. 
\end{rem}

\subsection{Riemannian structures.} Let us fix the group $G=\mathrm{O}(\R^n)$ of orthogonal matrices and $(M,\nabla,P)$ an infinitesimally homogeneous manifold with structure group $G$. It is well known that to have a $G$-structure in $M$ is equivalent to have a Riemannian metric in $M$. Since $-1\in G$, Lemma \ref{thm:lemasbasicos} implies that $T_0$ and $\mathfrak{J}_0$ vanish. Hence, $\nabla$ is the Levi-Civita connection associated to the Riemannian metric of $M$. Moreover, for the characteristic tensor $R_0$ it is clear that the procedure employed to solve the case $\mathrm{SO}(n)$ applies verbatim. 
In conclusion, the only suitable tensors for $R_0$ are the scalar multiples of the
following:
\begin{equation*} 
K:\R^n\times \R^n\times\R^n \ni (u,v,w)\longmapsto \langle v,w\rangle u -\langle u,w \rangle v \in \R^n.\end{equation*} 
Note that $R_0 = \lambda K$ is the curvature tensor of a Riemannian manifold of constant sectional curvature $\lambda$. That is, the infinitesimally homogeneous manifolds with structure group $G$ are the Riemannian manifolds endowed with the Levi-Civita connection and with constant sectional curvature.

\subsection{Volume structures.} Let us fix $G =  \mathrm{SL}(\R^n)$. Any $G$-structure $P\to M$ is determined by a volume form  $\omega$ in $M$, and vice versa. Let us consider $\nabla$ a linear connection in $TM$. The quotient $\mathfrak{gl}(T_xM)/\mathfrak{g}_x$ canonically identifies with $\mathbb R$ be means of the trace morphism. This, the inner torsion $\mathfrak J^P$ of the structure is a $1$-form in $M$. On the other hand, there is no linear non zero $G$-invariant $1$-form in $\mathbb R^n$. Hence, if $\mathfrak J^P$ is $G$-constant  then it vanishes. Equivalently we may say that the volume form $\omega$ is horizontal for the connection $\nabla$, that is, $\nabla\omega = 0$. Similarly, a linear algebra computation using generic parabolic maps in $\mathrm{SL}(\R^n)$ shows that there are not $G$-invariant tensors $T_0$ and $R_0$. Thus, if $(M,\nabla,P)$ is an infinitesimally homogeneous manifold, then $\nabla$ is a flat torsion-free linear connection. Summarizing, an infinitesimally homogeneous space with structure group $\mathrm{SL}(\R^n)$ is an affine manifold (in the sense of Thurston, \cite{HT}) endowed with a horizontal volume form.

\subsection{Constant rank distributions.} Let us fix $G=\mathrm{GL}(\R^n;\R^s)$ with $s<n$, the group of linear automorphisms of $\R^n$ that fix $\R^s$. It consist of block triangular matrices. Let $(M,\nabla,P)$ be an infinitesimally homogeneous manifold with structure group $G$. For this particular group, $G$-structures
$P\to M$ are in one-to-one correspondence with rank $s$ 
regular distributions of vector spaces in $M$. Let us denote by $L$ the distribution corresponding to $P$. The group
$G$ containts the center of $\mathrm{GL}(\R^n)$, hence by Lemma \ref{thm:lemasbasicos},
the characteristic tensors $\mathfrak J_0$, $T_0$ and $R_0$ vanish. We have that $(M,\nabla)$ is an affine manifold 
and $\nabla$ is a connection in $P$, or equivalently, the distribution $L$ is parallel with respect to $\nabla$, and hence it determines a foliation
$\mathcal F$ whose leaves are $s$-dimensional affine submanifolds $M$. Summarizing, the data of a infinitesimal homogeneous structure in $M$ with structure group $G$ consists of:
\begin{itemize}
\item[(a)] An affine manifold $(M,\nabla)$.
\item[(b)] A foliation $\mathcal F$ of $M$ by affine submanifolds of dimension $s$. 
\end{itemize}
Let us remind that a submanifold of $M$ is affine if $\nabla$ induces a torsion free and flat connection on it. Thus, the leaves of $\mathcal F$ are totally geodesic submanifolds of $M$. 

\subsection{Webs.} Let us fix $G \subset \rm{GL}(\R^n)$ the group of diagonal matrices, and let $(M,\nabla,P)$ be an infinitesimally homogeneous manifold with structure group $G$. By definition,  $G$-structure in $M$ is a web. As before, by Lemma \ref{thm:lemasbasicos}, the characteristic tensors $\mathfrak J_0$, $T_0$ and $R_0$ vanish. 
We have then that $(M,\nabla)$ is an affine manifold. The web $P$ is compatible with the connection
$\nabla$ and thus it is a web by geodesic lines. Summarizing, an infinitesimally homogeneous manifold whose structure group is that of diagonal matrices is an affine manifold endowed with a web (in the sense of \cite{AG}) of geodesic lines. 

\subsection{Almost hermitian structures.} Let us consider the group $G=\mathrm{U}(\R^{2n})$ of unitary matrices, and \linebreak $(M,\nabla,P)$ 
an infinitesimally homogeneous manifold with structure group $G$. It is well known (see propositions 3.1, 3.2 and 4.7 in \cite[vol. II]{KN}) that the $G$ structure $P$ determines an almost hermitian structure consisting in an almost complex structure $J$ and an hermitian metric $h$ in $M$. The bundle $P\to M$ is the bundle of unitary frames of the almost hermitian structure. Let us compute the 
characteristic tensors of $(M,\nabla,P)$. Let us proceed as in the case of Riemannian structures: since $-1\in G$, and Lemma \ref{thm:lemasbasicos}
implies that $T_0$ and $\mathfrak J_0$ vanish. It follows that $\nabla$ is the Levi-Civita connection of the hermitian metric $h$.
Let us recall the following result \cite[vol. II, Ch. IX, Theorem 4.8]{KN}:

\begin{teo}
Let $(M,J,h)$ be an almost hermitian manifold and $P\to M$ the bundle of unitary frames. Then $(M,J,h)$ is a K\"ahler manifold if and only if $P$ admits a torsion free connection (which is  necessarily unique).
\end{teo}

The vanishing of $\mathfrak J_0$ and $T_0$ implies that $\nabla$ is a torsion-free connection in $P$, thus $(M,J,h)$ is a K\"ahler manifold. Let us compute the curvature tensor $R_0$, let us consider the cuadrilinear map:
$${\bf K}\colon \mathbb C^n \times \mathbb C^n \times \mathbb C^n \times \mathbb C^n \to \mathbb R,$$
defined by the formula:
$${\bf K}(u_1,u_2,u_3,u_4) = $$
$$\frac{1}{4}\left( \langle u_1,u_3 \rangle \langle u_2,u_4 \rangle  - \langle u_1,u_4 \rangle  \langle u_2,u_3 \rangle + \langle u_1, iu_3 \rangle \langle u_2, iu_4 \rangle - \langle u_1,iu_4 \rangle \langle u_2,iu_3 \rangle + 2\langle u_1,iu_2 \rangle \langle u_3,iu_4 \rangle \right)$$

We have the following result analogous to Lemma \ref{curvature1}.

\begin{lema}\label{curvature2}
Let $R$ be a quadrilinear map $R\colon \mathbb R^n \times \mathbb R^n \times \mathbb R^n \times \mathbb R^n \longmapsto \mathbb R$ satisfying:
\begin{itemize}
\item[(a)] $R(u_1,u_2,u_3,u_4) = -R(u_1,u_2,u_3,u_4),$
\item[(b)] $R(u_1,u_2,u_3,u_4) = -R(u_1,u_2,u_4,u_3),$
\item[(c)] $R(u_1,u_2,u_3,u_4) + R(u_1,u_3,u_4,u_2) + R(u_1,u_4,u_2,u_3) = 0,$
\end{itemize}
for all $u_1,u_2,u_3,u_4\in \mathbb R^n$. If $R$ is ${\rm U}(n)$-invariant then $R$ is a multiple of the tensor ${\bf K}$.
\end{lema}

\begin{dem}
Let us consider $e_1,\ldots,e_n$ the canonical basis of $\mathbb C^n$ es complex vector space, so that $e_1,\ldots, e_n,$ $ie_1, \ldots, ie_n$ is a real basis. Let us take $\lambda = R(e_1,ie_1,e_1,ie_1)$. For vector $u\in \mathbb C^n$
there is an unitary transformation $g$ such that $g(u) = \|u\|e_1$ and thus:
$$R(u,iu,u,iu) = \langle u,u \rangle^2 R(e_1,ie_1,e_1,ie_1) = \lambda \langle u,u \rangle^2  = \lambda {\bf K}(u,iu,u,iu).$$
Finally, by Lemma \cite[vol II, Chapter IX, Proposition 7.1]{KN} we have that $R = \lambda {\bf K}$.
\end{dem}

Thus, we have $R_0 = \lambda{\bf K}$. This tensor $R_0$ is, by definition, the curvature tensor of a K\"ahler manifold of constant holomorphic curvature $\lambda$. Hence, $(M,\nabla,P)$ is a K\"ahler manifold of constant holomorphic curvature. Summarizing:

\begin{teo}
Let $(M,\nabla,P)$ be a manifold endowed with a connection $\nabla$ and an $U(n)$-structure $P$. It is an infinitesimally homogeneous manifold if and only if $P$ is the bundle of unitary frames of a constant holomorphic curvature K\"ahler structure in $M$ and $\nabla$ is its corresponding Levi-Civita connection. 
\end{teo}

\subsection{Product of Riemannian structures} Let us fix the group $\mathrm{O}(\R^{n_1})\times \mathrm{O}(\R^{n_2})\subset \GLnR$, where $n=n_1+n_2$. Since $-1\in G$, Lemma \ref{thm:lemasbasicos} implies that $T_0$ and $\mathfrak{J}_0$ vanish. Moreover, it follows from the procedure done for the orthogonal group that the only suitable tensors for the curvature are the linear combinations of $K$, $K_1$ and $K_2$ where for each $i=1,2$, $K_i$ is the version of $K$, defined in \eqref{eq:curvaturaconstante}, over $\R^{n_i}$.

On the other hand, since $\mathfrak{J}_0$ vanishes, it follows from the algebraic relation between $\mathfrak{J}_0$ and $R_0$ given in \cite[Theorem 1]{marin}
that a suitable candidate for the curvature tensor must be 
$(\mathfrak{o}(\R^{n_1})+\mathfrak{o}(\R^{n_2}))$-valued. Given a linear combination $aK + bK_1 + cK_2$, it is clear that it satisfies this condition only if $a$ vanishes. Thus,
the only candidates for the curvature are the 
linear combinations of $K_1$ and $K_2$.
Summarizing, an infinitesimally homogeneous manifold $(M,\nabla, P)$ with structural group $\mathrm{O}(\R^{n_1})\times \mathrm{O}(\R^{n_2})$ is locally described as the product $(M,g)=(M_1,g^1)\times (M_2,g^2)$, of two Riemannian manifolds $(M_1,g^1)$, $(M_2,g^2)$ each endowed with the Levi-Civita connection and with constant sectional curvature. Where $g$ denotes the metric structure given by the orthogonal sum of $g_1$ and $g_2$, more precisely, 
\[g_{(x_1,x_2)}\big((v_1,v_2),(w_1,w_2)\big)=g^1_{x_1}(v_1,w_1)+g^2_{x_2}(v_2,w_2),\]
for each $x_1\in M_1$, $x_2\in M_2$, $v_1,w_1\in T_{x_1}M_1$ and each
  $v_2,w_2\in T_{x_2}M_2$. Moreover, $\nabla$ is the Levi-Civita connection for $g$ and for each $x=(x_1,x_2) \in M$, $P_{(x_1,x_2)}$ is the set consisting of all linear isometries $p:\R^{n_1+n_2}\to T_{(x_1,x_2)}M$ such that
  \begin{gather*}
  p\big(\R^{n_1}\oplus\{0\}^{n_2}\big)=T_{x_1}M_1\oplus\{0\};\\ 
  p\big(\{0\}^{n_1}\oplus\R^{n_2}\big)=\{0\}\oplus
  T_{x_2}M_2.
  \end{gather*}

\section{Conclusions}

In this article we have discussed the infinitesimally homogeneous structures of affine manifolds with prescribed structure groups. We found that, for many examples we arrive to well known geometric structures. However, we used \emph{ad hoc} methods for each group. It would be interesting to consider some other examples, as symplectic groups, or develop a general procedure that allows to solve this problem for a more general family of matrix groups.

\section*{Acknowledgements}

We thank to the anonymous referees for their suggestions that allowed us to improve the scientific quality of this work. C. M. A. acknowledges the support of Universidad de Antioquia. D. B.-S. acknowledges the support of Universidad Nacional de Colombia through project ``M\'etodos algebraicos-geom\'etricos en ecuaciones diferenciales y temas relacionados 2015'' (Hermes 27984).


\begin{thebibliography}{9999}

\bibitem{AG} {\sc M. A. Akivis, V. V. Goldberg}, 
Differential geometry of webs. 
\emph{Handbook of differential geometry, 1} (2000), 1-152.


\bibitem{HT} {\sc M. Hirsch, W. Thurston}, 
\emph{Foliated bundles, invariant measures, and flat manifolds}, 
Ann. Math. (2) 101, (1975) 369-390.


\bibitem{KN} {\sc S. Kobayashi, K. Nomizu}, 
{\em Foundations of Differential Geometry},
vols.\ I,II New York, John Wiley \& Sons, Inc.\ (1963).

\bibitem{marin} {\sc Mar\'in C.}, \emph{An Algebraic Characterization of Affine Manifolds with $G$-structure Satisfying a Homogeneity Condition}, Revista Colombiana de Matem\'aticas \textbf{2} 2010, no. 2, 149-166.

\bibitem{tauskpaper} {\sc Piccione P., Tausk D.}, \emph{An Existence Theorem for G-structure Preserving Affine Immersions}, Indiana Univ. Math. J \textbf{57} 2008, no. 3, 1431-1465.

\bibitem{tauskbook} {\sc Piccione P., Tausk, D.}, 
{\em The theory of connections and $G$--structures: Applications to Affine and Isometric
Immersion}, XIV Escola de Geometria Diferencial, IMPA 2006.

\end{thebibliography}
\end{document}